\documentclass[bj]{imsart}
\startlocaldefs
\newtheorem{thm}{Theorem}[section]

\endlocaldefs

\startlocaldefs
\endlocaldefs

\begin{document}

\begin{frontmatter}

\title{Statistical generalized  derivative applied to the profile likelihood estimation in a mixture of semiparametric models}

\runtitle{Statistical generalized  derivative applied to a mixture of semiparametric models}

\begin{aug}
\author{\fnms{Yuichi} \snm{Hirose}\thanksref{a}\ead[label=e1]{Yuichi.Hirose@msor.vuw.ac.nz}}
\and
\author{\fnms{Ivy} \snm{Liu}\thanksref{a}\ead[label=e2]{i-ming.liu@vuw.ac.nz}}

\address[a]{ School of Mathematics, Statistics and Operations Research, Victoria University of Wellington,
 New~Zealand. 
\printead{e1,e2}
}

\runauthor{Yuichi Hirose}

\affiliation{Victoria University of Wellington}

\end{aug}

\begin{abstract}

There is a difficulty in finding an estimate of variance of the profile likelihood estimator in the joint model of longitudinal and survival data. 
We solve the difficulty by introducing the ``statistical generalized  derivative''.
The derivative is used to show the asymptotic normality of the estimator without assuming the second derivative of the density function in the model exists.

\end{abstract}

\begin{keyword}
\kwd{Efficiency}
\kwd{Efficient information bound}
\kwd{Efficient score}
\kwd{Implicitly defined function}
\kwd{Profile likelihood}
\kwd{Semi-parametric model}
\kwd{Joint model}
\kwd{EM algorithm}
\kwd{Mixture model}
\end{keyword}





\end{frontmatter}

\section{Introduction} 
This paper proposes a method to show asymptotic normality of a profile likelihood estimator in a mixture of semiparametric models with the EM-algorithm.
As an example we consider a joint model of ordinal responses and  the proportional hazards model with the finite mixture.
Through this example, we demonstrate to  solve the theoretical challenge in a joint model of survival and longitudinal data stated by \cite{Hsieh}:   
`` No distributional or asymptotic theory is available to date, and even the standard errors (SE), 
defined as the standard deviations of the parametric estimators, are difficult to obtain.''
The difficulty of the problem is to deal with an implicit function which is difficult to differentiate.
In the profile likelihood approach we profile out the baseline hazard function by plugging in an estimate of the hazard function to the likelihood function. 
This estimator of the hazard function is an implicit function in our problem.  

The core of our method is  an introduction of ``statistical generalised derivative'' (in Theorem 2.1).
Using this generalised derivative, in Theorem 2.2, we show asymptotic normality of estimator without differentiating the implicit function.   
In section 3, we apply our proposed method to the joint model.

Our approach gives an alternative to the methodologies  in \cite{Hirose11}, \cite{Hirose16} and \cite{MV00}, 
where an asymptotic normality of the profile likelihood estimator were studied.
Other related work is in \cite{Zeng05}. 
In this paper they showed asymptotic normality of the estimators through the joint maximization of the parameter of interest and the baseline hazard function.
This approach does not require to deal with the implicit function which encounter with the profile likelihood estimation.

\section{Mixture of semiparametric models and generalized statistical derivative}

 We consider a mixture of semiparametric models whose density is of the form
 \begin{eqnarray}\label{Mixture_dnsity}
 p(x;\theta,\eta,\pi)=\sum_{r=1}^R\pi_r p_r(x;\theta_r, \eta_r), 
 \end{eqnarray}
 where for each $r=1,\ldots,R$, $p_r(x;\theta_r, \eta_r)$ is a semiparametric model with a finite dimensional parameter $\theta_r \in \Theta_r \subset R^{m_r}$
 and an infinite dimensional parameter $\eta_r \in H_r$ where $H_r$ is a subset of Banach space $\mathcal{B}_r$, and $\pi_1,\ldots,\pi_R$ are mixture probabilities.
 We assume that $\pi_r>0$ for each $r$ and $\sum_{r=1}^R\pi_r =1$.
 We denote $\theta=(\theta_1,\ldots,\theta_R) \in \Theta=\Theta_1\times \cdots \times \Theta_R$, $\eta=(\eta_1,\ldots,\eta_R) \in H= H_1 \times \cdots \times H_R$ and $\pi=(\pi_1,\ldots,\pi_R)$.
 The true values of these parameters are denoted by
 $\theta_0=(\theta_{1,0},\ldots,\theta_{R,0})$, $\eta_0=(\eta_{1,0},\ldots,\eta_{R,0})$ and $\pi_0=(\pi_{1,0},\ldots,\pi_{R,0})$
Once we observe iid data $X_1,\ldots,X_n$ from the mixture model, the joint probability function of  the data ${\bf X}=(X_1,\ldots,X_n)$ is given by
\begin{eqnarray}\label{Lik_mix}
p({\bf X}; \theta,\eta,\pi)=\prod_{i=1}^n\sum_{r=1}^R\pi_r p_r(X_i;\theta_r, \eta_r). 
\end{eqnarray}
 We consider $\theta$ is the parameters of interest, and $\eta$ and $\pi$ are nuisance parameters.
This paper aims to establish large sample properties of the maximum likelihood estimator of $\theta$ using profile likelihood and the EM-algorithm (\cite{EM77}).

 To discuss the EM-algorithm, we further introduce notations (we use notations from \cite{Bishop06}).
 Let $Z_i=(Z_{i1},\ldots,Z_{iR})$ be group indicator variable for the subject $i$: for each $r$, $Z_{ir}=0$ or $=1$ with $P(Z_{ir}=1)=\pi_r$,  and $\sum_{r=1}^R Z_{ir}=1$.
 Let ${\bf Z}=(Z_1,\ldots,Z_n)$. The joint probability function of the complete data $({\bf X},{\bf Z})$ is
 \begin{eqnarray}\label{Lik_comp}
 p({\bf X},{\bf Z}; \theta,\eta,\pi)=\prod_{i=1}^n\prod_{r=1}^R[\pi_r p_r(X_i;\theta_r, \eta_r)]^{Z_{ir}}.  
 \end{eqnarray}

Then the EM-algorithm utilizes the identity
\begin{eqnarray}\label{EMeqn}
 \log p({\bf X};\theta,\eta,\pi) & = &  \sum_{{\bf Z}} q({\bf Z}) \log p({\bf X},{\bf Z}; \theta,\eta,\pi) 
 - \sum_{{\bf Z}} q({\bf Z}) \log p({\bf Z}|{\bf X}; \theta,\eta,\pi),
\end{eqnarray}
where $q({\bf Z})$ is any distribution of ${\bf Z}$ (\cite{EMbook}, Equation (3.3)).

In the E-step,
$$q({\bf Z})=p({\bf Z}|{\bf X}; \theta^{old},\eta^{old},\pi^{old}),$$
 then it is well known that the gradient for the $\log p({\bf X};\theta,\eta,\pi)$ coincides with 
 the one for $\sum_{{\bf Z}} q({\bf Z})\log p({\bf X},{\bf Z}; \theta,\eta,\pi)$  at $(\theta^{old},\eta^{old},\pi^{old})$. 
  In the M-step, maximize  the expectation of the complete data log likelihood function $\sum_{{\bf Z}} q({\bf Z})\log p({\bf X},{\bf Z}; \theta,\eta,\pi)$
  to obtain $(\theta^{new},\eta^{new},\pi^{new})$.
Then repeat E-step and M-step iteratively until we achieve the maximum. 

Under this procedure, 
the maximizer of the mixture log likelihood function $\log p({\bf X};\theta,\eta, \pi)$ 
with respect to $\theta$, $\eta$ and $\pi$ is the same as the ones for 
the expectation of the complete data log likelihood function $  \sum_{{\bf Z}} q({\bf Z}) \log p({\bf X},{\bf Z}; \theta,\eta,\pi)$ 
 (\cite{EMbook}, Section 3.4.1).

 
The EM-algorithm gives us the maximum likelihood estimator $\hat{\theta}$ of the mixture model.
However it does not give us the variance of the estimator.
In the following, we aim to establish asymptotic normality of the maximum likelihood estimator of $\theta$ using the profile likelihood estimation with the EM-algorithm.

\subsection{Generalized statistical derivative and asymptotic normality of the estimator}
From the complete data joint distribution (\ref{Lik_comp}), we can derive the conditional distribution $p({\bf Z}|{\bf X}; \theta,\eta,\pi)$:
\begin{eqnarray}\label{dist_Z|X}
 p({\bf Z}|{\bf X}; \theta,\eta,\pi) 
 & = & \frac{p({\bf X},{\bf Z};\theta,\eta,\pi)}{\sum_{\bf Z} p({\bf X},{\bf Z};\theta,\eta,\pi)}\nonumber \\
& = &\prod_{i=1}^n\prod_{r=1}^R \frac{[\pi_r p_r(X_i;\theta_r, \eta_r)]^{Z_{ir}}}{ \sum_{j=1}^R\pi_j p_j(X_i;\theta_j, \eta_j)}\nonumber \\
& = &\prod_{i=1}^n\prod_{r=1}^R \gamma_{r}(X_i;\theta,\eta)^{Z_{ir}}.
\end{eqnarray}
where 
\begin{eqnarray}\label{E(Z|X)}
\gamma_{r}(X_i;\theta,\eta)=\frac{\pi_r p_r(X_i;\theta_r, \eta_r)}{ \sum_{j=1}^R\pi_j p_j(X_i;\theta_j, \eta_j)}, \ \ r=1,\ldots,R. 
\end{eqnarray}

Again from (\ref{Lik_comp}), 
the expected complete data log-likelihood under $q({\bf Z})=p({\bf Z}|{\bf X}; \theta,\eta,\pi)$ is 
\begin{eqnarray}\label{complete_lik}
  \sum_{{\bf Z}} q({\bf Z}) \log p({\bf X},{\bf Z}| \theta,\eta,\pi)
 & = & \sum_{i=1}^n\sum_{r=1}^R \gamma_{r}(X_i;\theta,\eta)[ \log \pi_r + \log  p_r(X_i;\theta_r, \eta_r)].
 \end{eqnarray}

With the expected complete data log-likelihood  (\ref{complete_lik}),
the method of Lagrange multiplier  can be applied to get the MLE $\hat{\pi}_k$ of $\pi_r$:
\begin{eqnarray}\label{hat_pi}
\hat{\pi}_k(\theta,\eta) =\frac{\sum_{i=1}^n \gamma_{r}(X_i;\theta,\eta)}{n}, \ \ r=1,\ldots,R. 
\end{eqnarray}
We require that, as $n \rightarrow \infty$,
$$\hat{\pi}_r(\theta_0,\eta_0) \stackrel{P}{\rightarrow} \pi_{r,0}$$
where $(\theta_0,\eta_0)$ are the true value of $(\theta,\eta)$ and $\pi_{r,0}$, $r=1,\ldots,R$, are the true mixture probabilities.

\bigskip
{\bf The efficient score function and information matrix in the mixture model:} 
 The score function for $\theta$ and score operator for $\eta$ in the mixture model given in (\ref{Mixture_dnsity}) are, respectively,
 \begin{eqnarray}\label{score1}
\dot{\ell}(x;\theta,\eta) 
& = & \frac{\partial}{\partial \theta}\log \left(\sum_{r=1}^R\pi_r p_r(x;\theta_r, \eta_r)\right)
 =  \sum_{r=1}^R \gamma_{r}(x;\theta,\eta)\frac{\partial}{\partial \theta} \log p_r(x;\theta_r, \eta_r),
 \end{eqnarray}
 and
 \begin{eqnarray}\label{score2}
 B(x;\theta,\eta)  & = & d_{\eta}\log \left(\sum_{r=1}^R\pi_r p_r(x;\theta_r, \eta_r)\right)
 =  \sum_{r=1}^R \gamma_{r}(x;\theta,\eta)d_{\eta} \log p_r(x;\theta_r, \eta_r)
 \end{eqnarray}
where $\gamma_{r}(x;\theta,\eta)$ is given in (\ref{E(Z|X)}) with $X_i$ replaced with $x$. 
The notation $d_{\eta}$ is the Hadamard derivative operator with respect to the parameter $\eta$.

Let $\theta_0,\eta_0$ be the true values of $\theta,\eta$ and denote 
$\dot{\ell}_0(x)=\dot{\ell}(x;\theta_0,\eta_0)$ and $B_0(x)=B(x;\theta_0,\eta_0)$.
Then, it follows from the standard theory (\cite{Vaart98}, page 374) that the efficient score function $\tilde{\ell}_0$ and the efficient information matrix $\tilde{I}_0$ in the semiparametric mixture model are given by
\begin{eqnarray}\label{eff_score}
\tilde{\ell}_0(x)=(I-B_0(B_0^*B_0)^{-1}B_0^*)\dot{\ell}_0(x), 
\end{eqnarray}
 and
 \begin{eqnarray}\label{eff_information}
\tilde{I}_0=E[\tilde{\ell}_0\tilde{\ell}_0^T].  
 \end{eqnarray}

 {\bf Note:} 
 Equations (\ref{score1}) and (\ref{score2}) show that the score functions in the semiparametric mixture model (\ref{Mixture_dnsity}) coincide with the ones for 
 the expected complete data likelihood (\ref{complete_lik}).

 \bigskip
{\bf The score function for the profile likelihood:} 
In the estimation of $(\theta,\eta)$ we use  the profile likelihood approach:
we obtain a function $(\theta,F) \rightarrow \hat{\eta}_{\theta,F} =(\hat{\eta}_{1,\theta,F},\ldots,\hat{\eta}_{R,\theta,F} ) $ 
whose values are in the space of the parameter $\eta=(\eta_1,\ldots,\eta_R)$.

Define the score functions for the profile likelihood in the model
\begin{eqnarray}\label{Profile_score}
 \phi(x;\theta,F)= \frac{\partial}{\partial \theta}\log \left(\sum_{r=1}^R\pi_r p_r(x;\theta_r, \hat{\eta}_{r,\theta,F})\right)
 = \sum_{r=1}^R  \gamma_{r}(x;\theta,\hat{\eta}_{\theta,F}) \frac{\partial}{\partial \theta} \log  p_r(x;\theta_r, \hat{\eta}_{r,\theta,F})
\end{eqnarray}
and
\begin{eqnarray}\label{Profile_score2}
 \psi(x;\theta,F)= d_{F}\log \left(\sum_{r=1}^R\pi_r p_r(x;\theta_r, \hat{\eta}_{r,\theta,F})\right)
 = \sum_{r=1}^R  \gamma_{r}(x;\theta,\hat{\eta}_{\theta,F}) d_F \log  p_r(x;\theta_r, \hat{\eta}_{r,\theta,F}),
\end{eqnarray}
We require that $\eta_0=\hat{\eta}_{\theta_0,F_0}=(\hat{\eta}_{1,\theta_0,F_0},\ldots,\hat{\eta}_{R,\theta_0,F_0})$ and
the condition (R2) below assumes $\phi(x;\theta_0,F_0)$ is the efficient score function $\tilde{\ell}_0(x)$ in the model where
$\theta_0$, $\eta_0$ and $F_0$ are the true values of the parameters $\theta$, $\eta$ and cdf $F$.

\bigskip
{\bf Assumptions:} We list assumptions used for Theorem \ref{theorem1} and Theorem \ref{theorem2} given below.

On the set of cdf functions ${\cal F}$, we use the sup-norm, i.e.
for $F, F_0 \in {\cal F}$, $$\|F-F_0\|=\sup_x|F(x)-F_0(x)|.$$
For $\rho>0$, let
$${\cal C}_{\rho}=\{F \in {\cal F}:\|F-F_0\| < \rho\}.$$

We assume that:
\begin{enumerate}

\item[(R1)] 
For each $(\theta,F) \in \Theta \times {\cal F}$, the log-profile-likelihood function for an observation $x$
\begin{eqnarray}\label{log-density}
\log p(x;\theta,F)=\log \left(\sum_{r=1}^R\pi_r p_r(x;\theta_r, \hat{\eta}_{r,\theta,F})\right)
\end{eqnarray}
is continuously differentiable with respect to $\theta  = (\theta_1,\ldots, \theta_R)$ and Hadamard differentiable with respect to $F$ for all $x$.
Derivatives are respectively denoted by  $\phi(x;\theta,F)= \frac{\partial}{\partial \theta}\log p(x;\theta,F)$ 
and $\psi(x;\theta,F)= d_F\log p(x;\theta,F)$ and they are given  in (\ref{Profile_score}) and  (\ref{Profile_score2}). 

\item[(R2)]  
We denote $\hat{\eta}_{\theta,F}=(\hat{\eta}_{1,\theta,F},\ldots,\hat{\eta}_{R,\theta,F})$.
We assume $\hat{\eta}_{\theta,F}$ satisfies $\hat{\eta}_{\theta_0,F_0}=\eta_0=(\eta_{1,0},\ldots,\eta_{R,0})$ and the function
$$\tilde{\ell}_0(x):=\phi(x;\theta_0,F_0)$$
is the efficient score function.
Further, we assume the cube-root-n consistency: if $\hat{\theta}_n$ is the MLE of $\theta_0$,
$n^{1/3}(F_n-F_0)=O_P(1)$ and $n^{1/3}(\hat{\eta}_{\hat{\theta}_n,F_n}-\eta_0)=O_P(1)$.


\item[(R3)] The efficient information matrix $\tilde{I}_0=E[\tilde{\ell}_0\tilde{\ell}_0^{T}]
=E[\phi\phi^{T}(X;\theta_0,F_0)]$ is invertible.

\item[(R4)] 
The score function $\phi(x;\theta,F)$ defined in (\ref{Profile_score}) takes the form
$$\phi(x;\theta,F)=\tilde{\phi}(x;\theta,F,\hat{\eta}_{\theta,F}),$$
where, by assumption (R2), the efficient score function is given by 
$$\tilde{\ell}_0(x)=\phi(x;\theta_0,F_0)=\tilde{\phi}(x;\theta_0,F_0,\eta_0).$$ 

We assume that there exists a $\rho>0$ and neighborhoods $\Theta$ and $H$ of $\theta_0$  and  $\eta_0$, respectively, 
 such that ${\cal C}_{\rho}$ and $H$ are Donsker and 
the class of functions $\{\tilde{\phi}(x;\theta,F, \eta):\ (\theta,F, \eta) \in \Theta \times {\cal C}_{\rho} \times H \}$ has a square integrable envelope function and
it is Lipschitz in the parameters $(\theta,F, \eta)$:
\begin{eqnarray}\label{tilde_phi_Lipschitz}
\|\tilde{\phi}(x;\theta',F', \eta')-\tilde{\phi}(x;\theta,F, \eta)\| \leq M'(x) (\|\theta' -\theta\|+\|F'-F\|+\|\eta'-\eta\|) 
\end{eqnarray}
where $M'(x)$ is a $P_0$-square integrable function.
Moreover, for $\theta,\theta' \in \Theta$ and $F,F' \in {\cal C}_{\rho}$,
\begin{eqnarray}\label{p_Lipschitz}
\left|\frac{p(x;\theta',F')-p(x;\theta,F)}{p(x;\theta,F)}\right| & \leq  & M(x) (\|\theta'-\theta\| +\|F'-F\|)
\end{eqnarray}
where $M(x)$ is a $P_0$-square integrable function.

\end{enumerate}

\bigskip
{\bf Main result: statistical generalized derivative and asymptotic linearity of the estimator.}
To calculate the second derivative of the score function $\phi(x;\theta,F)$ given in (\ref{Profile_score}), 
we use the idea similar to the derivative of generalized functions (\cite{Kolmogorov}).
Let $\varphi \rightarrow (f,\varphi)=\int_{-\infty}^{\infty}f(x) \varphi(x) dx$ be a generalized function, where $\varphi$ vanishes outside of some interval.
Then if $f$ and $\varphi$ are differentiable with derivative $f'$ and $\varphi'$, then by integration by parts,
$$(f',\varphi)=\int_{-\infty}^{\infty}f'(x) \varphi(x) dx=-\int_{-\infty}^{\infty}f(x) \varphi'(x) dx=-(f,\varphi').$$
We define the derivative $(f',\varphi)$ of the generalized function $\varphi \rightarrow (f,\varphi)$ by $-(f,\varphi')$.
This definition is valid even if $f$ is not differentiable, provided $\varphi$ is differentiable.

A similar idea can be applied in our problem.
Suppose the density for the profile likelihood $p(x;\theta,F)$ given in (\ref{log-density}) 
is twice differentiable with respect to $\theta$, then by differentiating 
\begin{eqnarray*}
 \int  
 \left\{\frac{\partial}{\partial \theta} \log p(x;\theta,F)\right\}
 p(x;\theta,F)dx=0,
\end{eqnarray*}
with respect to $\theta$ at $(\theta,F)=(\theta_0,F_0)$, 
we get equivalent expressions for the efficient information matrix in terms of the score function $\phi(x;\theta_0,F_0)$:
\begin{eqnarray}\label{eff_information2}
\tilde{I}_0=E [\phi \phi^T(X;\theta_0,F_0)]= - E\left[\frac{\partial}{\partial \theta^T}\phi(X;\theta_0,F_0) \right].
\end{eqnarray}
From this equation we are motivated to define the expected derivative of the score function 
$- E\left[\frac{\partial}{\partial \theta^T}\phi(X;\theta_0,F_0) \right]$ by $E [\phi \phi^T(X;\theta_0,F_0)]$.
In the following theorem, we show that the definition is valid even when the derivative of the score function 
$\frac{\partial}{\partial \theta^T}\phi(x;\theta,F)$ does not exist.

\begin{thm}\label{theorem1}
Suppose (R1) and (R4).
Let $p(x;\theta,F)=\sum_{r=1}^R\pi_r p_r(x;\theta_r, \hat{\eta}_{r,\theta,F})$, 
 $\phi(x;\theta,F)= \frac{\partial}{\partial \theta}\log p(x;\theta,F)$,
 and $\psi(x;\theta,F)= d_F\log p(x;\theta,F)$ as defined in (\ref{log-density}), (\ref{Profile_score}) and (\ref{Profile_score2}), respectively.
Let $\theta_t$ and $F_t$ be a smooth paths through $\theta_0$ and $F_0$ at $t=0$ such that
 the limits of $t^{-1}(\theta_t - \theta_0)$ and $t^{-1}(F_t - F_0)$ exist as $t \rightarrow 0$.
Then,  as $ t \rightarrow 0$, we have that
\begin{eqnarray}\label{limit1}
   & & E \left[t^{-1} \{\phi(X;\theta_t,F_0)-\phi(X;\theta_0,F_0)\}\right] \nonumber \\
   & & \hspace{0.5cm}  = -E\left[\phi(X;\theta_0,F_0) \phi^T(X;\theta_0,F_0)\right]  \{t^{-1}(\theta_t-\theta_0)\}+o(1),
\end{eqnarray}
and
\begin{eqnarray}\label{limit2}
  & &   E \left[t^{-1} \{\phi(X;\theta_t,F_t)-\phi(X;\theta_t,F_0)\}\right] \nonumber \\
   &  & \hspace{0.5cm} =  - E[ \phi(X;\theta_0,F_0) \psi(X;\theta_0,F_0)] \{t^{-1}(F_{t}-F_0)\}\nonumber \\
& & +o(1)+O\{(t^{-1}\|\theta_t-\theta_0\|+t^{-1}\|F_{t}-F_0\|) (\|F_t-F_0\|+\|\hat{\eta}_{\theta_t,F_t}-\hat{\eta}_{\theta_t,F_0}\|)\}.
\end{eqnarray}
   
\end{thm}

{\bf Note.}
Note that even when the derivative $\frac{\partial}{\partial \theta} \phi(x;\theta,F)$ does not exist the equation (\ref{equality}) in the proof holds.  
Together with the derivative $\frac{\partial}{\partial \theta}p(x;\theta,F)$ exists implies 
that the derivative of the map $\theta \rightarrow  E \left[\phi(x;\theta,F)\right] $ exists and it is given by (\ref{limit1}).
We may call the derivative the statistical generalized derivative.
A similar comment for (\ref{limit2}) holds.

{\bf Proof.}
We assumed the limits of $t^{-1}(\theta_t - \theta_0)$ and $t^{-1}(F_t - F_0)$ exist as $t \rightarrow 0$. 
By the differentiability of $p(x;\theta,F)$ with respect to $\theta$ and $F$, 
at each $x$ with $p(x;\theta_0,F_0)>0$ we have, as $t \rightarrow 0$,
\begin{eqnarray}\label{deriv1}
\frac{t^{-1} \{ p(x;\theta_t,F_{0}) - p(x;\theta_0,F_0)}{p(x;\theta_0,F_0)} =\phi(x;\theta_0,F_0) \{t^{-1}(\theta_{t}-\theta_0)\}+o(1),  
\end{eqnarray}
and
\begin{eqnarray}\label{deriv2}
\frac{t^{-1} \{ p(x;\theta_t,F_{t}) - p(x;\theta_t,F_0)}{p(x;\theta_0,F_0)} =\psi(x;\theta_0,F_0) \{t^{-1}(F_{t}-F_0)\}+o(1).  
\end{eqnarray}

We prove (\ref{limit1}).
For each $t$, the equality
\begin{eqnarray*}
 0 & = &  t^{-1}\left\{\int \phi(x;\theta_t,F_0) p(x;\theta_t,F_0)dx- \int \phi(x;\theta_0,F_0) p(x;\theta_0,F_0)dx\right\}\\ 
  & = &  \int t^{-1} \{\phi(x;\theta_t,F_0)-\phi(x;\theta_0,F_0)\} p(x;\theta_0,F_0)dx\\ 
  &  &  +  \int \phi(x;\theta_t,F_0) t^{-1} \{ p(x;\theta_t,F_0) - p(x;\theta_0,F_0)\}dx
\end{eqnarray*}
holds, where we understood the integral is taken over the set $\{x: p(x;\theta_0,F_0)>0\}$.
It follows that, for each $t$, we have that
\begin{eqnarray}\label{equality}
& &     \int t^{-1} \{\phi(x;\theta_t,F_0)-\phi(x;\theta_0,F_0)\} p(x;\theta_0,F_0)dx \nonumber \\
&  & \hspace{1cm}   =  -   \int \phi(x;\theta_t,F_0) t^{-1} \{ p(x;\theta_t,F_0) - p(x;\theta_0,F_0)\}dx.
\end{eqnarray}

By Appendix 1 (a), the right hand side of (\ref{equality}) is, as $t \rightarrow 0$,
\begin{eqnarray*}
& & -   \int \phi(x;\theta_t,F_0) t^{-1} \{ p(x;\theta_t,F_0) - p(x;\theta_0,F_0)\}dx\\
& = & -   \int \phi(x;\theta_t,F_0) \frac{t^{-1} \{ p(x;\theta_t,F_0) - p(x;\theta_0,F_0)\}}{p(x;\theta_0,F_0)}p(x;\theta_0,F_0)dx\\
 & = &  - \int \phi(x;\theta_0,F_0) \phi^T(x;\theta_0,F_0)p(x;\theta_0,F_0)dx \left\{ t^{-1}(\theta_t-\theta_0)\right\} +o(1).
\end{eqnarray*}

It follows that, we have (\ref{limit1}):
\begin{eqnarray*}
&  & \int t^{-1} \{\phi(x;\theta_t,F_0)-\phi(x;\theta_0,F_0)\} p(x;\theta_0,F_0)dx\\
&  & \hspace{1cm} = - \int \phi(x;\theta_0,F_0) \phi^T(x;\theta_0,F_0)p(x;\theta_0,F_0)dx \{t^{-1}(\theta_t-\theta_0)\} +o(1).
\end{eqnarray*}

\bigskip
Now we prove (\ref{limit2}).
Similar to the beginning of the proof of  (\ref{limit1}), for each $t$, the following equation holds: 
\begin{eqnarray}\label{eqn2}
 &  &  \int t^{-1} \{\phi(x;\theta_t,F_{t})-\phi(x;\theta_t,F_0)\} p(x;\theta_t,F_{t})dx \nonumber \\
    &   & \hspace{0.5 cm}= - \int \phi(x;\theta_t,F_0) t^{-1} \{ p(x;\theta_t,F_{t}) - p(x;\theta_t,F_0)\}dx.
\end{eqnarray}

By Appendix 1 (b), the left hand side of (\ref{eqn2}) is, as $t \rightarrow 0$,
\begin{eqnarray}\label{equality3}
& &  \left\| \int t^{-1} \{\phi(x;\theta_t,F_{t})-\phi(x;\theta_t,F_0)\} p(x;\theta_t,F_{t})dx 
-  \int t^{-1} \{\phi(x;\theta_t,F_{t})-\phi(x;\theta_t,F_0)\} p(x;\theta_0,F_{0})dx\right\| \nonumber \\
& = & \left\| \int  \{\phi(x;\theta_t,F_{t})-\phi(x;\theta_t,F_0)\} \frac{t^{-1}\{p(x;\theta_t,F_{t}) -  p(x;\theta_0,F_{0})\}}{p(x;\theta_0,F_{0})} p(x;\theta_0,F_{0})dx \right\| \nonumber \\
&  =  &   O\{t^{-1}(\|\theta_t-\theta_0\|+\|F_{t}-F_0\|) (\|F_t-F_0\|+\|\hat{\eta}_{\theta_t,F_t}-\hat{\eta}_{\theta_t,F_0}\|)\}. 
\end{eqnarray}

Using  (\ref{deriv2}), the similar proof of Appendix 1 (a) can show  that the integral in the right hand side of the equation (\ref{eqn2}) is 
\begin{eqnarray}\label{equality4}
& &   \int \phi(x;\theta_t,F_0) t^{-1} \{ p(x;\theta_t,F_{t}) - p(x;\theta_t,F_0)\}dx\nonumber \\
& = &  \int \phi(x;\theta_0,F_0) \psi(x;\theta_0,F_0) t^{-1}(F_t-F_0)   p(x;\theta_0,F_0)dx +o(1). 
\end{eqnarray}

By combining (\ref{equality3}) and (\ref{equality4}), the equality (\ref{eqn2}) is equivalent to
\begin{eqnarray*}
& &   \int t^{-1} \{\phi(x;\theta_t,F_{t})-\phi(x;\theta_t,F_0)\}  p(x;\theta_0,F_0)dx \\
& = & - \int \phi(x;\theta_0,F_0) \psi(x;\theta_0,F_0)p(x;\theta_0,F_0)dx \{t^{-1}(F_{t}-F_0)\}\\
& & +o(1) +O\{t^{-1}(\|\theta_t-\theta_0\|+\|F_{t}-F_0\|) (\|F_t-F_0\|+\|\hat{\eta}_{\theta_t,F_t}-\hat{\eta}_{\theta_t,F_0}\|)\}.
\end{eqnarray*}
The (\ref{limit2}) follows from this.


\bigskip
Using the result in Theorem \ref{theorem1}, we show the following result:
\begin{thm}\label{theorem2}
Suppose the set of assumptions $(R1)-(R4)$ holds. 
Then a consistent solution $\hat{\theta}_n$ to the estimating equation 
\begin{eqnarray}\label{ProfileEstEqn}
\sum_{i=1}^{n} \phi(X_{i};\hat{\theta}_n,F_n)=0
\end{eqnarray}
 is an asymptotically linear estimator for $\theta_0$ :
$$\sqrt{n}(\hat{\theta}_n -\theta_0)=\frac{1}{\sqrt{n}} \sum_{i=1}^{n} \tilde{I}_0^{-1} \tilde{\ell}_0(X_{i}) + o_P(1).$$
Hence we have that
$$\sqrt{n}(\hat{\theta}_n -\theta_0) \stackrel{d}{\longrightarrow} N\left(0,\tilde{I}_0^{-1} \right) \ \ \textrm{as } n\rightarrow \infty.$$
\end{thm}

{\bf Proof}

In (R4) we assumed 
${\cal C}_{\rho}$ and $H$ are Donsker and 
the function $\tilde{\phi}(x;\theta,F, \eta)$
is Lipschitz in the parameters $(\theta,F, \eta)$ with a $P_0$-square integrable function $M'(x)$ given in (\ref{tilde_phi_Lipschitz}).
By Corollary 2.10.13 in \cite{VW96}, the class 
$\{\tilde{\phi}(x;\theta,F, \eta):\ (\theta,F, \eta) \in \Theta \times {\cal C}_{\rho} \times H \}$ 
is Donsker.

By Lemma 19.24 in \cite{Vaart98} together with the dominated convergence theorem, it implies
\begin{eqnarray}\label{Donsker}
 \frac{1}{\sqrt{n}} \sum_{i=1}^{n} \{\phi(X_{i};\hat{\theta}_n,F_n)-\phi(X_{i};\theta_0,F_0)\}
= \sqrt{n}E\{\phi(X;\hat{\theta}_n,F_n)-\phi(X;\theta_0,F_0)\}+o_P(1).
\end{eqnarray}

From  (\ref{limit1}) it follows that
\begin{eqnarray} \label{Proof_expand_1}
  \sqrt{n}E\{\phi(X;\hat{\theta}_n,F_0)-\phi(X;\theta_0,F_0)\} 
& = &   - \tilde{I}_0 \sqrt{n}(\hat{\theta}_n-\theta_0)+o_p(1),
\end{eqnarray}
where $\tilde{I}_0=E[ \tilde{\ell}_0 \tilde{\ell}_0^T]=E\{\phi(X;\theta_0,F_0)\phi^T(X;\theta_0,F_0)\}$.

Using (\ref{limit2}),
\begin{eqnarray}\label{Proof_expand_2}
   &  & \sqrt{n}E \{\phi(X;\hat{\theta}_n,F_n)-\phi(X;\hat{\theta}_n,F_0)\} \nonumber \\
   &    = &   - E[ \phi(X;\theta_0,F_0) \psi(X;\theta_0,F_0)] \{\sqrt{n}(F_{n}-F_0)\}  \nonumber \\
& &   +O\{\sqrt{n}(\|\hat{\theta}_n-\theta_0\|+\|F_{n}-F_0\|) (\|F_n-F_0\|+\|\hat{\eta}_{\hat{\theta}_n,F_n}-\eta_0\|)\} \nonumber\\
& & + o(1+\|\hat{\theta}_n-\theta_0\|+\|F_{n}-F_0\| +\|\hat{\eta}_{\hat{\theta}_n,F_n}-\eta_0\|) \nonumber \\
& =  &    o_P(1+ \sqrt{n}(\hat{\theta}_n-\theta_0)),
\end{eqnarray}
where we used:
\begin{enumerate}
\item Since $\psi(x;\theta_0,F_0)$ is in the nuisance tangent space and $\phi(x;\theta_0,F_0)$ is the efficient score function, we have
\begin{eqnarray}\label{eqn_efficient}
E[\phi(x;\theta_0,F_0) \psi(x;\theta_0,F_0)] =0. 
\end{eqnarray}
\item  We assumed $(\hat{\theta}_n-\theta_0)=o_P(1)$, $n^{1/3}(F_n-F_0)=O_P(1)$ and $n^{1/3}(\hat{\eta}_{\hat{\theta}_n,F_n}-\eta_0)=O_P(1)$, 
it follows that  
\begin{eqnarray*}
O\{\sqrt{n}(\|\hat{\theta}_n-\theta_0\|+\|F_{n}-F_0\|) (\|F_n-F_0\|+\|\hat{\eta}_{\hat{\theta}_n,F_n}-\eta_0\|)\}  & = & o_P(1+ \sqrt{n}(\hat{\theta}_n-\theta_0)) \\
\textrm{ and }  \hspace{1cm} o(1+\|F_{n}-F_0\| +\|\hat{\eta}_{\hat{\theta}_n,F_n}-\eta_0\|)& = & o_P(1). 
\end{eqnarray*}

\end{enumerate}

Using (\ref{Proof_expand_1}) and  (\ref{Proof_expand_2}),  the right hand side of (\ref{Donsker}) is
\begin{eqnarray}\label{Stat_Deriv}
& &  \sqrt{n}E\{\phi(X;\hat{\theta}_n,F_n)-\phi(X;\theta_0,F_0)\} \nonumber\\
& = &   \sqrt{n}E\{\phi(X;\hat{\theta}_n,F_0)-\phi(X;\theta_0,F_0)\}+ \sqrt{n}E\{\phi(X;\hat{\theta}_n,F_n)-\phi(X;\hat{\theta}_n,F_0)\} \nonumber \\
& = &   - \tilde{I}_0 \sqrt{n}(\hat{\theta}_n-\theta_0)+o_p\{1+\sqrt{n}(\hat{\theta}_n-\theta_0)\}.
\end{eqnarray}
Finally,  
(\ref{Donsker}) together with (\ref{Stat_Deriv}) and  $\frac{1}{\sqrt{n}} \sum_{i=1}^{n} \phi(X_{i};\hat{\theta}_n,F_n)=0$ imply that
\begin{eqnarray*}
\sqrt{n}(\hat{\theta}_n-\theta_0)   = \frac{1}{\sqrt{n}} \sum_{i=1}^{n} \tilde{I}_0^{-1}\phi(X_{i};\theta_0,F_0)
+o_P(1).
\end{eqnarray*}

\section{Joint mixture model of survival and longitudinal ordered data}

In this section, we apply the theorem 2.1 and 2.2 to the example of  ``the joint model of ordinal responses and  the proportional hazards  with the finite mixture'' 
which is studied in  \cite{Preedalikit}.
We demonstrate that how our method can solve the difficulty in profile likelihood estimation in  the joint model.  

The maximum likelihood estimation in the joint model has been studied by many authors, among others we name few, \cite{Wulfsohn97}, \cite{Song02} and \cite{Hsieh}.
For more complete review of the joint models please see \cite{Tsiatis04} and \cite{Rizopoulos12}.

\bigskip
\noindent{\bf Ordinal Response Models:}
Let $Y_{ijm}$ be the ordered categorical response from $1$ (poor) to $L$ (excellent) on item (or question) $j$ for subject $i$ at the $m^{th}$ protocol-specified time point,
where $i=1,2,\ldots,n$, $j=1,2,\ldots,J$ and $m=1,2,\ldots,M$. 
In total, there are $J$ items in the questionnaire related to patients quality of life, collected at times $t_1,t_2,\ldots,t_M$. 
Given that subject $i$ belongs to group $r$, an ordered stereotype model can be written as
\begin{eqnarray*}
 \label{semi:2}
\log \left[\frac{P(Y_{ijm}= \ell \ | \ \theta_r)}{P(Y_{ijm}= 1 \ | \ \theta_r)}\right] = a_{\ell} + \phi_\ell(b_j+\theta_r),\,\,\, r=1,\ldots, R,
\end{eqnarray*}
where $a_{\ell}$ is a response level intercept parameter with
$\ell=2,\ldots, L$, $b_j$ is an item effect, and $\theta_r$ is
associated with the discrete latent variable, with $a_{1}=0$, $b_1 =
0$, $\phi_1=0$ and $\theta_1 = 0$.  The parameter $\theta_r$ can be
referred to as a group effect of the quality of life for patients in
group $r$. However, the group memberships are unknown. The
$\{\phi_\ell\}$ parameters can be regarded as unknown scores for the
outcome categories. Because
$\phi_\ell(b_j+\theta_r)=(A\phi_\ell((b_j+\theta_r)/A))$ for any
constant $A\ne 0$, for
identifiability, we need to impose monotone scores on
$\{\phi_\ell\}$ to treat $Y_{ijm}$ as ordinal. Therefore, the model has the constraint $0= \phi_1
\leq \phi_2 \leq \ldots\leq \phi_L = 1$. 
The ordinal response part of likelihood function for the $i$th subject is
\begin{eqnarray}
\label{semi:3}
P (Y_i \,| \, \theta_r , \alpha )  \, 
&=& \, \prod_{m=1}^{M_i}\prod_{j=1}^{J}\prod_{\ell=1}^{L}\bigg(\frac{\exp(a_{\ell}+\phi_\ell (b_j +\theta_r))}{1+\sum_{k=2}^{L}\exp(a_{k}+ \phi_k( b_j +\theta_r))}\bigg)^{Y_{ijm\ell}}
\end{eqnarray}
where $\alpha= (a,b,\phi)$.
Each follow-up time point may have a different number of observations because some patient responses are missing.

\bigskip
\noindent{\bf The Cox Proportional Hazards Model:}
We consider the Cox proportional hazards model for the
survival part in the joint model. 
Let $X$ be a
time-independent covariate.
 The
hazard function for the failure time $T_i$ of the $i^{th}$ subject is of the
form
\begin{eqnarray} \label{sur:1}
 \lambda(t | X_i, \theta_r, \delta) 
	& = & \lambda_0(t)\exp(\theta_r \delta_0 + X_i\delta_1) 
\end{eqnarray}
where $\lambda_0(t)$ is the baseline hazard function.
The latent variable $\theta_r$ is linked with the ordinal response model and
$\delta=(\delta_0,\delta_1)$ are coefficients. 

For the estimation of the baseline hazard function $\lambda_0(t)$, we use
the method of nonparametric maximum likelihood described in
\cite[section 4.3]{b16}. Let $\lambda_i$ be the hazard at time $t_i$, where
$t_1 < t_2 < \ldots < t_n$ are the ordered observed times. 
Assume that the hazard is zero between adjacent times so that the survival time is discrete. 
The corresponding cumulative hazard function $\displaystyle
\Lambda_0(t_i) = \sum_{p \leq i} \lambda_p$ is a step function with jumps at
the failure time $t_i$. Then the survival part likelihood function of subject $i$ is
\begin{eqnarray}
\label{cox_dist}
   P\big(T_i,d_i \,| \, \lambda  , \theta_r , \delta \big)
    & =  & \, \big(\lambda_i\exp(\theta_r \delta_0 + X_i\delta_1) \big)^{d_i} \times \exp\Big(- \sum_{p \leq i} \lambda_p \exp(\theta_r \delta_0 + X_i\delta_1)\Big),
\end{eqnarray}
where the $d_{i}$ is an indicator of censorship for
individual $i$:
if we observe failure time, then $d_{i}=1$, otherwise
$d_{i}=0$.

\bigskip
\noindent{\bf The Full Likelihood Function:}
The joint likelihood function is
obtained by combining the probability function from ordinal response model (\ref{semi:3}), and
the proportional hazards model (\ref{cox_dist}), by
assuming the two models are independent given latent discrete random
variables.

Let $\pi_r$ be the unknown
probability ($r=1,\ldots,R$) that a subject lies in group $r$, and
$(\Theta,\lambda)=((\theta,\alpha,\delta),\lambda)$ be all the unknown parameters of the joint model. The
mixture model likelihood function is 
\begin{eqnarray}
\label{full_likelihood}
L( \Theta,\lambda | Y, T,D) 
= \prod_{i=1}^{n} \Bigg(\sum_{r=1}^{R} P\big (Y_i \,| \, \theta_r , \alpha \big) P\big (T_i,d_i \,| \, \lambda  , \theta_r , \delta \big) \pi_r \Bigg).
\end{eqnarray}

Let $Z_{ir}$ be the group indicator, where $Z_{ir}$ = 1 if the $i^{th}$ individual was from the $r^{th}$ group and 0 otherwise. The complete data likelihood can be written as
\begin{eqnarray}
\label{comp}
L( \Theta,\lambda | Y, T,d,Z) 
= \prod_{i=1}^{n}\prod_{r=1}^{R} \Big(\, P\big (Y_i \,| \, \theta_r , \alpha \big) 
P\big (T_i,d_i \,| \, \lambda  , \theta_r , \delta \big) \pi_r\Big)^{Z_{ir}}. 
\end{eqnarray}

The expected complete data log likelihood under $q({\bf Z})=P({\bf Z}| Y, T,d)$ is
\begin{eqnarray}\label{Ecomp}
& & \sum_{\bf Z} q({\bf Z})\log L( \Theta, \lambda | Y, T,d,Z)\nonumber \\
	& = & \sum_{i=1}^{n}\sum_{r=1}^{R}  \gamma(Z_{ir}) 
	\left\{\log \pi_r  + \log P\big (Y_i \,| \, \theta_r , \alpha \big) 
+ \log P\big (T_i,d_i \,| \, \lambda  , \theta_r , \delta \big) \right\}
\end{eqnarray}
where $\gamma(Z_{ir})$, $P\big (Y_i \,| \, \theta_r , \alpha \big)$ and $P\big (T_i,d_i \,| \, \lambda  , \theta_r , \delta \big)$
are defined in equations (\ref{semi:8}),  (\ref{semi:3}) and (\ref{cox_dist}) respectively.

To estimate all parameters and the baseline hazards simultaneously, we
combine the EM algorithm and the method of nonparametric maximum
likelihood.

\subsection{Estimation procedure: profile likelihood with EM algorithm}

\noindent{\bf Baseline Hazard Estimation:}
Before starting the EM-step, we profile out the baseline hazard
function $\lambda_0(t)$.  The survival part of equation (\ref{Ecomp}) can be separately maximized with respect to $\lambda$:
\begin{eqnarray}
\label{semi:9}
& & \sum_{i=1}^{n}\sum_{r=1}^{R}  \gamma(Z_{ir})  \log P\big (T_i,d_i \,| \, \lambda  , \theta_r , \delta \big) \nonumber \\
 & &  \hspace{1cm}=  \sum_{i=1}^{n}\sum_{r=1}^{R}  \gamma(Z_{ir})  \Big\{
d_i(\log\lambda_i + \theta_r \delta_0 + X_i\delta_1)   - \sum_{p \leq i} \lambda_p \exp(\theta_r \delta_0 + X_i\delta_1) \Big\}.
\end{eqnarray}


By solving $\frac{\partial}{\partial \lambda_l}\sum_{i=1}^{n}\sum_{r=1}^{R}  \gamma(Z_{ir})  \log P\big (T_i,d_i \,| \, \lambda  , \theta_r , \delta \big)=0 $, $l=1,\ldots,n$,
we find the maximizer $\widehat{\lambda}_l$ of (\ref{semi:9}) by holding $(\theta , \delta)$ fixed, and it is given by 
\begin{eqnarray}
\label{hat_lambda}
 \widehat{\lambda}_l(\theta,\delta) &=  \frac{d_i}{\sum_{p\geq i}\sum_{r=1}^{R}  \gamma(Z_{pr}) \exp(\theta_r \delta_0 + X_p\delta_1)}.
\end{eqnarray}

Denote $\widehat{\lambda}(\theta,\delta) =(\widehat{\lambda}_1(\theta,\delta),\ldots, \widehat{\lambda}_n(\theta,\delta))$.

\noindent{\bf The E-step:}
In the E-step, we use the current parameter estimates $\Theta = (\theta, \alpha , \delta) $ to find the expected values of $Z_{ir}$:
\begin{eqnarray}\label{semi:8}
 \gamma(Z_{ir})=E\big(Z_{ir}  | \, Y_i , T_i , d_i \big) 
	&  =  & \frac{\pi_r \,P\big(Y_i  \,|\, \theta_r, \alpha \big) P\big (T_i,d_i \,| \,  \widehat{\lambda}(\theta,\delta) , \theta_r , \delta \big)}
	{\sum_{g=1}^{R} \pi_g \,P\big(Y_i  \,|\, \theta_g , \alpha\big) P\big (T_i,d_i \,| \,\widehat{\lambda}(\theta,\delta) , \theta_g , \delta \big)}.
\end{eqnarray}

\noindent{\bf The M-step:}
In the M-step, we maximize equation (\ref{Ecomp})   with respect to $\pi_r$ and $\Theta = (\theta, \alpha , \delta) $. Due to the
fact that there is no relationship between $\pi_r$ and $\Theta$,
they can be estimated separately.
\begin{enumerate}
\item
Calculate the estimates of $\pi_r$
\begin{eqnarray*}
\widehat{\pi_r} & = &  \frac{\sum_{i=1}^{n} \gamma(Z_{ir})}{n} . 
\end{eqnarray*}

\item We maximize the second and third parts of equation (\ref{Ecomp}) (with $ \widehat{\lambda}(\theta,\delta)$ in the place of $\lambda$) 
  \begin{eqnarray}\label{Ecomp2}
	 \sum_{i=1}^{n}\sum_{r=1}^{R}  \gamma(Z_{ir} ) 
	 \left\{   \log P\big (Y_i \,| \, \theta_r , \alpha \big) 
	 + \log P\big (T_i,d_i \,| \, \widehat{\lambda}(\theta,\delta) , \theta_r , \delta \big) \right\}
\end{eqnarray}
   with respect to $\Theta =
  (\theta, \alpha , \delta)$ to obtain $\widehat{\Theta}$. 
  
\end{enumerate}
The estimated parameters from the M-step are returned into the E-step until the value of $\widehat{\Theta}$ converges.

\subsection{Asymptotic normality of the MLE $\widehat{\Theta}$ and its asymptotic variance}

From (\ref{hat_lambda}), an estimator of the cumulative hazard function 
in the counting process notation is 
$$\widehat{\Lambda}(t)
=\int_0^t  \frac{  \sum_{i=1}^n d N_i(u)}{\sum_{i=1}^n Y_i(u) \sum_{r=1}^{R}  \gamma(Z_{ir}) \exp(\theta_r \delta_0 + X_i\delta_1)}  $$
where $N_i(u)=1_{\{T_i \leq u, d_i=1\}}$ and $Y_i(u)=1_{\{T_i \geq u\}}$.

Let us denote $E_{F_n}f=\int f dF_n$. Then the above $\widehat{\Lambda}(t)$ can be written as
\begin{eqnarray}\label{hat_Lambda}
\widehat{\Lambda}(t;\Theta,F_n)
=\int_0^t  \frac{ E_{F_n} dN(u) }{E_{F_n}  Y(u) \sum_{r=1}^{R}  \gamma(Z_{r}) \exp(\theta_r \delta_0 + X\delta_1)}   
\end{eqnarray}
where $N(u)=1_{\{T \leq u, d=1\}}$, $Y(u)=1_{\{T \geq u\}}$ and
similarly $\gamma(Z_{r})$ is defined.

Equation (\ref{Ecomp2}) gives the profile likelihood function for $\Theta =
  (\theta, \alpha , \delta)$. 
The log-profile likelihood function for one observation is
\begin{eqnarray}\label{LogP_joint}
\log P(Y_i , T_i , d_i|\Theta,F_n)=  \sum_{r=1}^{R}  \gamma(Z_{ir})\left\{ \log P\big (Y_i \,| \, \theta_r , \alpha \big)
+ \log P\big (T_i,d_i \,| \, \widehat{\Lambda}(\Theta,F_n) , \theta_r , \delta \big)\right\} 
\end{eqnarray}
where
\begin{eqnarray}\label{LogP_1}
& & \sum_{r=1}^{R}  \gamma(Z_{ir})\log P (Y_i \,| \, \theta_r , \alpha ) \\
&  & \hspace{0.5cm} = \sum_{r=1}^{R}  \sum_{m=1}^{M_i}\sum_{j=1}^{J}\sum_{\ell=1}^{L} \gamma(Z_{ir}) Y_{ijm\ell} \left\{a_{\ell}+\phi_\ell (b_j +\theta_r)
-\log \left(1+\sum_{k=2}^{L}\exp(a_{k}+ \phi_k( b_j +\theta_r))\right)\right\},\nonumber
\end{eqnarray}
and
\begin{eqnarray}\label{LogP_2}
& &\sum_{r=1}^{R} \gamma(Z_{ir}) \log P\big (T_i,d_i \,| \, \widehat{\Lambda}(\Theta,F_n)  , \theta_r , \delta \big)\nonumber\\ 
 & = &     
   \sum_{r=1}^R\gamma(Z_r)\left\{d_i  \left(\log \frac{ E_{F_n} dN(T_i) }{E_{F_n}  Y(T_i) \sum_{r'=1}^{R}  \gamma(Z_{r'}) \exp(\theta_{r'} \delta_0 + X\delta_1)}+ \theta_r \delta_0 + X_i\delta_1 \right)\right. \nonumber\\
    & &\left. -  \exp(\theta_r \delta_0 + X_i\delta_1) \int_0^{T_i}  \frac{ E_{F_n} dN(u) }{E_{F_n}  Y(u) \sum_{r'=1}^{R}  \gamma(Z_{r'}) \exp(\theta_{r'} \delta_0 + X\delta_1)} \right\}.
\end{eqnarray}
In the above log-likelihood we set $a_1=b_1=\phi_1=\theta_1=0$.

\bigskip
{\bf Score functions}

The score functions for the profile likelihood are
\begin{eqnarray}\label{score_joint_profile}
& & \phi(Y_i , T_i , d_i|\Theta,F_n)= \phi_O(Y_i |\Theta)+ \phi_S(T_i , d_i|\Theta,F_n) \nonumber\\  
& & \hspace{0.5cm}  =    \sum_{r=1}^{R}  \gamma(Z_{ir})  \frac{\partial}{\partial \Theta} \log P\big (Y_i \,| \, \theta_r , \alpha \big) 
	 + \sum_{r=1}^{R}  \gamma(Z_{ir})  \frac{\partial}{\partial \Theta}\log P\big (T_i,d_i \,| \, \widehat{\Lambda}(\Theta,F_n) , \theta_r , \delta \big),\nonumber\\
& & \psi(Y_i , T_i , d_i|\Theta,F_n)=    
	  \sum_{r=1}^{R}  \gamma(Z_{ir})  d_F \log P\big (T_i,d_i \,| \, \widehat{\Lambda}(\Theta,F_n) , \theta_r , \delta \big).
	 \end{eqnarray}
	 Here all derivatives are calculated treating $\gamma(Z_{ir})$ as constant.
We call
 $\phi_O$ is the score function for the ordinal response model and $\phi_S$ is the one for the survival model.

\begin{thm}\label{theorem3} (The efficient score function)
We drop subscript $i$ in equation (\ref{score_joint_profile}).
We have the followings: at the true value of $(\Theta, F)$,
\begin{enumerate}
 \item $ \widehat{\Lambda}(t;\Theta,F) =\Lambda(t)$, the true cumulative hazard function,  and;
 \item the score function $\phi(Y , T , d|\Theta,F)$ defined in (\ref{score_joint_profile})  is the efficient score function in the model.
\end{enumerate}
\end{thm}

The proof of Theorem 3.1 
is given in Appendix 2.

 \subsubsection{Checking conditions}
 We check conditions (R1)-(R4) in Section 2.1 so that Theorem 2.1 and 2.2 can be used to get the large sample distribution of the estimator $ \widehat{\Theta}_n$:
 $$\sqrt{n}(\widehat{\Theta}_n- \Theta) \sim N(0,\tilde{I}^{-1}),$$
 where $\tilde{I}=E(\phi \phi^T)$ is the efficient information with $\phi$ is defined in (\ref{score_joint_profile}).
  
Since the ordinal response data part is a parametric model, we mainly discuss for the survival part of the model.
The survival part of the profile log -likelihood function for a one observation is given in (\ref{LogP_2}).

To express the survival part of the score function $\phi_S (T,d|\Theta,F)$ in the form given in condition (R4), we introduce a few notations. 

Let
\begin{eqnarray}\label{gamma_Z}
  \gamma(Z_{r}|\Theta, \Lambda)
	&  =  & \frac{\pi_r \,P\big(Y  \,|\, \theta_r, \alpha \big) P\big (T,d \,| \, \Lambda , \theta_r , \delta \big)}
	{\sum_{g=1}^{R} \pi_g \,P\big(Y  \,|\, \theta_g , \alpha \big) P\big (T,d \,| \, \Lambda  , \theta_g , \delta \big)}.
\end{eqnarray}
The function $ \gamma(Z_{r}|\Theta, \Lambda)$ is differentiable with respect to $\Theta$ and  $\Lambda$.
Then the function $\gamma(Z_r)$ in (\ref{LogP_2}) can be expressed as
$$\gamma(Z_r)=\gamma(Z_{r}|\Theta, \widehat{\Lambda}(\Theta,F)).$$

Let
\begin{eqnarray*}
M_0(t|\Theta,F,\Lambda) & = & E_F  Y(u) \sum_{r=1}^{R}  \gamma(Z_{r}|\Theta,\Lambda) \exp(\theta_r \delta_0 + X\delta_1)\\
M_1(t|\Theta,F,\Lambda) & = & E_F  Y(u) \sum_{r=1}^{R}  \gamma(Z_{r}|\Theta,\Lambda) \left(\begin{array}{c}\delta_0 \\ \theta_r \\ X\end{array}\right)   \exp(\theta_r \delta_0 + X\delta_1).
\end{eqnarray*}
Then the score function for the survival part $\phi_S (T,d|\Theta,F )$  is
\begin{eqnarray}\label{Score_survial}
& &\tilde{\phi}_S (T,d|\Theta,F ,\widehat{\Lambda}(\Theta,F)) \nonumber \\
  &  = &  \sum_{r=1}^R  \gamma(Z_{r}|\Theta, \widehat{\Lambda}(\Theta,F)) \frac{\partial}{\partial \Theta}\log P\big (T,d \,| \, \widehat{\Lambda}(\Theta,F)  , \theta_r , \delta \big) \nonumber \\
   & = &  
   \sum_{r=1}^R \gamma(Z_{r}|\Theta, \widehat{\Lambda}(\Theta,F))\left\{d  \left[ \left(\begin{array}{c}\delta_0 \\ \theta_r \\ X\end{array}\right) - \frac{M_1(T|\Theta,F, \widehat{\Lambda}(\Theta,F))}{M_0(T|\Theta,F, \widehat{\Lambda}(\Theta,F))}
   \right]\right.  \\
    & &\left. 
      +   \exp(\theta_r \delta_0 + X\delta_1) \int_0^T  \frac{ E_F dN(u)
      }{M_0(u|\Theta,F, \widehat{\Lambda}(\Theta,F))}  
       \left[\left(\begin{array}{c}\delta_0 \\ \theta_r \\ X\end{array}\right) - \frac{M_1(u|\Theta,F, \widehat{\Lambda}(\Theta,F))}{M_0(u|\Theta,F, \widehat{\Lambda}(\Theta,F))}\right] 
    \right\}.\nonumber
\end{eqnarray}

We will check condition (R4)  using the function defined by
\begin{eqnarray}\label{Score_surval2}
& & \tilde{\phi}_S (T,d|\Theta,F,\Lambda) \nonumber \\
   & = &  
   \sum_{r=1}^R \gamma(Z_{r}|\Theta, \Lambda)
   \left\{d  \left[ \left(\begin{array}{c}\delta_0 \\ \theta_r \\ X\end{array}\right) 
   - \frac{M_1(T|\Theta,F, \Lambda)}{M_0(T|\Theta,F, \Lambda)}
   \right]\right.  \\
    & &\left. 
      +   \exp(\theta_r \delta_0 + X\delta_1) \int_0^T  \frac{ E_F dN(u)
      }{M_0(u|\Theta,F, \Lambda)}  
       \left[\left(\begin{array}{c}\delta_0 \\ \theta_r \\ X\end{array}\right) - \frac{M_1(u|\Theta,F, \Lambda)}{M_0(u|\Theta,F, \Lambda)}\right] 
    \right\}.\nonumber
\end{eqnarray}

{\bf Condition (R1):}
We calculated the survival part score function  $\phi_S(T,d|\Theta,F)=\tilde{\phi}_S (T,d|\Theta,F ,\widehat{\Lambda}(\Theta,F)) $ in (\ref{Score_survial}).
The ordinal response data part is a parametric model, it is differentiable with respect to the parameter $\Theta$ (we omit the calculation).

We calculate the score function $\psi(T , d|\Theta,F)= \sum_{r=1}^{R}  \gamma(Z_{ir})  d_F \log P\big (T,d \,| \, \widehat{\Lambda}(\Theta,F) , \theta_r , \delta \big)$.
For an integrable function $h$ with the same domain as the cdfs $F$,
\begin{eqnarray*}
\psi(Y , T , d|\Theta,F)h
  & = &  \sum_{r=1}^R\gamma(Z_r)d_F \log P\big (T,d \,| \, \widehat{\Lambda}(\Theta,F)  , \theta_r , \delta \big)h\\
   & = &  
   \sum_{r=1}^R\gamma(Z_r)\left\{d  \left(\frac{E_hdN(T)  }{ E_F dN(T)  }  - \frac{ E_hY(T) \sum_{r=1}^{R}  \gamma(Z_{r}) \exp(\theta_r \delta_0 + X\delta_1)}{ E_{F}  Y(T) \sum_{r=1}^{R}  \gamma(Z_{r}) \exp(\theta_r \delta_0 + X\delta_1)}
    \right)\right.\\
    & & -  \exp(\theta_r \delta_0 + X\delta_1) \int_0^T  \frac{ E_{h} dN(u)}{E_F  Y(u) \sum_{r=1}^{R}  \gamma(Z_{r}) \exp(\theta_r \delta_0 + X\delta_1)}  \\
    & & \left. +  \exp(\theta_r \delta_0 + X\delta_1)\int_0^T \frac{E_{F} dN(u) E_h  Y(u) \sum_{r=1}^{R}  \gamma(Z_{r}) \exp(\theta_r \delta_0 + X\delta_1)}{(E_F  Y(u) \sum_{r=1}^{R}  \gamma(Z_{r}) \exp(\theta_r \delta_0 + X\delta_1))^2}\right\}.
\end{eqnarray*}



\bigskip
{\bf Condition (R2):}
We assume cube root $n$ consistency of the empirical cdf and an estimate of the baseline cumulative hazard function: $\|F_n-F_0\|=O_p(n^{-1/3})$ and
$\|\widehat{\Lambda}(\widehat{\Theta}_n,F_n)-\Lambda\|=O_p(n^{-1/3})$ (where the both norms are the sup norm).

In Theorem 3.1 we verified the rest of conditions in (R2).

\bigskip
{\bf Condition (R3):}
We outline verification of condition (R3)

We treat $y_{ijm\ell}, T_i, d_i, X_i$ as random and the rest as constants. It requires some inspections to see that there is no linear combination of derivatives listed below is constant:
$\frac{\partial }{\partial a_k}\sum_{r=1}^{R}  \gamma(Z_{ir})\log P (Y_i \,| \, \theta_r , \alpha )$  ,
$\frac{\partial }{\partial \phi_{k}}\sum_{r=1}^{R}  \gamma(Z_{ir}) \log P (Y_i \,| \, \theta_r , \alpha )$,\\
$\frac{\partial }{\partial b_j} \sum_{r=1}^{R}  \gamma(Z_{ir})\log P (Y_i \,| \, \theta_r , \alpha ) $,
$\frac{\partial  }{\partial \theta_r} \sum_{r=1}^{R}  \gamma(Z_{ir})\log P(Y_i \,| \, \theta_r , \alpha )$,\\
$\frac{\partial}{\partial \delta_p}\sum_{r=1}^{R}  \gamma(Z_{ir}) \log P (T_i,d_i \,| \, \widehat{\lambda}(\theta,\delta)  , \theta_r , \delta)$
($k=2,\ldots,L; j=2,\ldots,J; r=2,\ldots,R; p=0,1$).
To ensure to this happens we put $a_1=b_1=\phi_1=\theta_1=0$.
It follows that the score function (\ref{score_joint_profile}) has an invertible variance-covariance matrix (cf. Theorem 1.4 in \cite{Seber:Lee:2003}).

\bigskip
{\bf Condition (R4):}
The score function $\tilde{\phi}_S (T,d|\Theta,F,\Lambda)$ given in 
(\ref{Score_surval2}) is differentiable with respect to the parameters $(\Theta,F,\Lambda)$ we assume the derivatives are bounded by square integrable envelope functions.
It follows that the score function is Lipschitz in parameters.

We also assume that the density in the model given in (\ref{LogP_joint}) satisfy (\ref{p_Lipschitz}).

\section{Discussion}

The proposed ``statistical generalized  derivative'' in Theorem 2.1 is applied for the score function in the model.
In this approach we do not require differentiability of the score function to show the asymptotic normality of the profile likelihood estimator in the model (in Theorem 2.2). 
However, it still requires the differentiability of the density function in the model.
In our joint model example in the section 3, 
the efficient score function was calculated without differentiating the implicit function (Theorem 3.1, the calculation is in Appendix 1).
In the example, we established asymptotic normality of the estimator without differentiating the implicit function.

There may be some examples that require to differentiate an implicit function in the calculation of the efficient score function. 
If this is the case the approach in \cite{Hirose16} may be applicable.
To demonstrate this approach in the joint model example in the paper, in Appendix 3,
we proved the differentiability of the implicit  function $\widehat{\Lambda}(t;\Theta,F)$ given in (\ref{hat_Lambda}).

Once we have the efficient score function of the model under consideration, we can apply the Theorem 2.1 and 2.2 in the paper 
to show the asymptotic normality of the profile likelihood estimator in the problem.


\section*{Appendix 1: Dominated convergence and bound for integral in the proof of Theorem 2.1}
Let $\theta_t$ and $F_t$ be a smooth paths through $\theta_0$ and $F_0$ at $t=0$ such that
 the limits of $t^{-1}(\theta_t - \theta_0)$ and $t^{-1}(F_t - F_0)$ exist as $t \rightarrow 0$.
Under the assumptions of Theorem 2.1, (R1)--(R4), we show the followings:
\begin{enumerate}
 \item[(a)]
  \begin{eqnarray*}
& &  \int  \phi(x;\theta_{t},F_0)  \frac{t^{-1}\{p(x;\theta_{t},F_0)-p(x;\theta_0,F_0) \}}{p(x;\theta_0,F_0)} p(x;\theta_0,F_0)dx \\
& =  &  \int   \phi(x;\theta_0,F_0) \phi^T(x;\theta_0,F_0)\{t^{-1}(\theta_{t}-\theta_0)\} p(x;\theta_0,F_0)dx +o(1) \textrm{ as }  t \rightarrow 0.
  \end{eqnarray*}
 
 \item[(b)] As $t \rightarrow 0$,
  \begin{eqnarray*}
& &   \left\| \int  \{\phi(x;\theta_t,F_{t})-\phi(x;\theta_t,F_0)\} \frac{t^{-1}\{p(x;\theta_t,F_{t}) -  p(x;\theta_0,F_{0})\}}{p(x;\theta_0,F_{0})} p(x;\theta_0,F_{0})dx\right\|\\
 & = &  O\{(t^{-1}\|\theta_t-\theta_0\|+ t^{-1}\|F_t-F_0\|)(\|F_t-F_0\|+\|\hat{\eta}_{\theta_t,F_t}- \hat{\eta}_{\theta_t,F_0}\|) \}.
\end{eqnarray*}

 \end{enumerate}

{\bf Proof of (a).} 
 Let $a \in R^p$ be a fixed nonzero vector. 
 Since $p(x;\theta,F)$ is differentiable with respect to $\theta$,  for each $x$ with $p(x;\theta_0,F_0)>0$, 
 \begin{eqnarray*}
G(x;t,a) & := & a^T\phi(x;\theta_t,F_0)  \frac{t^{-1}\{p(x;\theta_t,F_0)-p(x;\theta_0,F_0) \}}{p(x;\theta_0,F_0)}\\
& = & a^T \phi(x;\theta_0,F_0) \phi^T(x;\theta_0,F_0)\{t^{-1}(\theta_t-\theta_0)\} +o(1) \textrm{ as } t \rightarrow 0.
 \end{eqnarray*}

 By assumption (R4), there is a $P_0$-square integrable function $M(x)$ such that
 \begin{eqnarray*}
  \|\phi(x;\theta,F)\| & \leq &  M(x)\\
  \bigg| \frac{t^{-1}\{p(x;\theta_t,F_0)-p(x;\theta_0,F_0) \}}{p(x;\theta_0,F_0)} \bigg| & \leq &  M(x) \| t^{-1}(\theta_t-\theta_0)\|
 \end{eqnarray*}

Then
$$|G(x;t,a)| \leq M^2(x) \|a\| \| t^{-1}(\theta_t-\theta_0)\|$$
 and
 \begin{eqnarray*}
  \int  a^T\phi(x;\theta_t,F_0)  \frac{t^{-1}\{p(x;\theta_t,F_0)-p(x;\theta_0,F_0) \}}{p(x;\theta_0,F_0)} p(x;\theta_0,F_0)dx =
   \int  G(x;t,a) p(x;\theta_0,F_0)dx, 
 \end{eqnarray*}
where we understood the integral is over the set $\{x: p(x;\theta_0,F_0)>0 \}$.

Let $t_n$ be a sequence such that $t_n \rightarrow 0^+$.
Since  $ M^2(x) \|a\|  \| t^{-1}(\theta_t-\theta_0)\| \pm  G(x;t,a)$ are nonnegative function, by the Fatou's lemma,
\begin{eqnarray*}
 &  &\left(\int  M^2(x) p(x;\theta_0,F_0)dx\right) \|a\|\liminf_n \| t_n^{-1}(\theta_{t_n}-\theta_0)\|+\int  \liminf_n G(x; t_n,a)  p(x;\theta_0,F_0)dx\\
 & = &\int  \liminf_n \left[M^2(x)\|a\| \| t_n^{-1}(\theta_{t_n}-\theta_0)\|+G(x; t_n,a) \right] p(x;\theta_0,F_0)dx\\
 & \leq & \liminf_n \int   \left[M^2(x)\|a\| \| t_n^{-1}(\theta_{t_n}-\theta_0)\|+ G(x; t_n,a) \right] p(x;\theta_0,F_0)dx\\
 & = &   \left(\int   M^2(x)  p(x;\theta_0,F_0)dx\right) \|a\|  \liminf_n\| t_n^{-1}(\theta_{t_n}-\theta_0)\|
 +  \liminf_n  \int  G(x;t_n,a)  p(x;\theta_0,F_0)dx.
\end{eqnarray*}
 Hence we have
 $$\int  \liminf_n G(x;t_n,a)  p(x;\theta_0,F_0)dx \leq  \liminf_n  \int  G(x; t_n,a)  p(x;\theta_0,F_0)dx.$$

 Similarly,  by the Fatou's lemma,
\begin{eqnarray*}
& & \left(\int   M^2(x) p(x;\theta_0,F_0)dx\right)\|a\| \liminf_n\| t_n^{-1}(\theta_{t_n}-\theta_0)\|- \int  \limsup_n  G(x;t_n,a)  p(x;\theta_0,F_0)dx\\
 & = &\int  \liminf_n \left(M^2(x)\|a\| \| t_n^{-1}(\theta_{t_n}-\theta_0)\|- G(x;t_n,a) \right) p(x;\theta_0,F_0)dx\\
 & \leq & \liminf_n \int   \left(M^2(x)\|a\| \| t_n^{-1}(\theta_{t_n}-\theta_0)\|-G(x; t_n,a) \right) p(x;\theta_0,F_0)dx\\
 & \leq &  \left(\int   M^2(x)  p(x;\theta_0,F_0)dx\right)\|a\|  \liminf_n\| t_n^{-1}(\theta_{t_n}-\theta_0)\| 
 -  \limsup_n  \int  G(x; t_n,a)  p(x;\theta_0,F_0)dx.
\end{eqnarray*}
 From this it follows that
 $$  \limsup_n  \int  G(x; t_n,a)  p(x;\theta_0,F_0)dx
 \leq 
 \int  \limsup_n G(x; t_n,a)  p(x;\theta_0,F_0)dx.$$

 Combine all then we have
 \begin{eqnarray*}
 & & \int  \liminf_n G(x; t_n,a)  p(x;\theta_0,F_0)dx \leq  \liminf_n  \int  G(x;t_n,a)  p(x;\theta_0,F_0)dx \\
& \leq &  \limsup_n  \int  G(x; t_n,a)  p(x;\theta_0,F_0)dx
 \leq 
 \int  \limsup_n G(x; t_n,a)  p(x;\theta_0,F_0)dx.
 \end{eqnarray*}
Since 
$$ \liminf_n G(x; t_n,a)=\limsup_n G(x; t_n,a)$$
 we get
 \begin{eqnarray*}
  \lim_n \int  G(x; t_n,a) p(x;\theta_0,F_0)dx
&  = &   \int  \lim_n G(x; t_n,a) p(x;\theta_0,F_0)dx.
 \end{eqnarray*}
 
Equivalently, we have 
 \begin{eqnarray*}
& &  \int  a^T\phi(x;\theta_{t},F_0)  \frac{t^{-1}\{p(x;\theta_{t},F_0)-p(x;\theta_0,F_0) \}}{p(x;\theta_0,F_0)} p(x;\theta_0,F_0)dx \\
& =  &  \int  a^T \phi(x;\theta_0,F_0) \phi^T(x;\theta_0,F_0)\{t^{-1}(\theta_{t}-\theta_0)\} p(x;\theta_0,F_0)dx +o(1)
 \end{eqnarray*}
 as $t \rightarrow 0$.
 Since the vector $a$ is arbitrary nonzero vector, we have shown (a).

 {\bf Proof of (b).} 
 
 Using  (\ref{tilde_phi_Lipschitz}) in (R4) with a $P_0$-square integrable function $M'(x)$, we have 
\begin{eqnarray}\label{phi_dif}
\|\phi(x;\theta_t,F_{t})-\phi(x;\theta_t,F_0)\| & = &  \|\tilde{\phi}(x;\theta_t,F_{t},\hat{\eta}_{\theta_t,F_t})-\tilde{\phi}(x;\theta_t,F_{0},\hat{\eta}_{\theta_t,F_0})\|\nonumber \\
& \leq & M'(x) (\|F_t-F_0\|+\|\hat{\eta}_{\theta_t,F_t}-\hat{\eta}_{\theta_t,F_0}\|).
\end{eqnarray}

By (\ref{p_Lipschitz}) in (R4), there is a $P_0$-square integrable function $M(x)$ such that
 \begin{eqnarray*}
  \bigg| \frac{t^{-1}\{p(x;\theta_t,F_t)-p(x;\theta_0,F_0) \}}{p(x;\theta_0,F_0)} \bigg| & \leq &  M(x)( \| t^{-1}(\theta_t-\theta_0)\|+ \|t^{-1}(F_t-F_0)\|).
 \end{eqnarray*}

 Using these 
 we get the result 
 \begin{eqnarray*}
& &   \left\| \int  \{\phi(x;\theta_t,F_{t})-\phi(x;\theta_t,F_0)\} \frac{t^{-1}\{p(x;\theta_t,F_{t}) -  p(x;\theta_0,F_{0})\}}{p(x;\theta_0,F_{0})} p(x;\theta_0,F_{0})dx\right\|\\
 & \leq & \left\| \int  M'(x) M(x) p(x;\theta_0,F_{0})dx\right\|\\
 & & \times (\|F_t-F_0\|+\|\hat{\eta}_{\theta_t,F_t}- \hat{\eta}_{\theta_t,F_0}\|)( \| t^{-1}(\theta_t-\theta_0)\|+ \|t^{-1}(F_t-F_0)\|)\\
 & = &  O\{(t^{-1}\|\theta_t-\theta_0\|+ t^{-1}\|F_t-F_0\|)(\|F_t-F_0\|+\|\hat{\eta}_{\theta_t,F_t}- \hat{\eta}_{\theta_t,F_0}\|) \}.
 \end{eqnarray*}

\section*{Appendix 2: Proof of Theorem 3.1 (The Efficient score function)}

{\bf Proof.}
From (\ref{hat_Lambda}), replacing $F_n$ by $F$, we have
$$\widehat{\Lambda}(t;\Theta,F)
=\int_0^t  \frac{ E dN(u) }{E  Y(u) \sum_{r=1}^{R}  \gamma(Z_{r}) \exp(\theta_r \delta_0 + X\delta_1)}  $$
where 
$E$ is the expectation with respect to the true distribution $F$.
Since, at the true value of the parameters $(\Theta,F,\Lambda)$,
   \begin{eqnarray}\label{EdN}
     E dN(u) & = & E [ Y(u) \sum_{r=1}^{R}  \gamma(Z_{r}) \exp(\theta_r \delta_0 + X\delta_1)]d\Lambda(u),
    \end{eqnarray}
we have that $ \widehat{\Lambda}(t;\Theta,F) =\Lambda(t)$.

	 The score function $\phi(Y,T,d|\Theta,F)=\phi_O(Y|\Theta)+\phi_S (T,d|\Theta,F)$ in (\ref{score_joint_profile}) has two parts: 
the score function for the ordinal response model $\phi_O(Y|\Theta)$
and the score function for the survival model
$\phi_S (T,d|\Theta,F)$.
Since the score function for the ordinal response model does not involve the parameter $\Lambda$, we will only work on the survival part of score function.

We treat the part $\gamma(Z_r)$ as constant in terms of the parameters.

Let 
\begin{eqnarray}\label{M1M0}
M_{1}(t)& = & E\sum_{r=1}^R\gamma(Z_r)\left(\begin{array}{c}\delta_0 \\ \theta_r \\ X \end{array}\right)  \exp(\theta_r \delta_0 + X\delta_1) I(t \leq T) \nonumber \\ 
M_{0}(t)& = & E\sum_{r=1}^R\gamma(Z_r)  \exp(\theta_r \delta_0 + X\delta_1) I(t \leq T) 
\end{eqnarray}

Then the score function in the survival part of the model at the true value of parameters $\Theta$ and $F$ is
\begin{eqnarray}\label{eff_score_survival}
  & &\phi_S (T,d|\Theta,F )= \sum_{r=1}^R\gamma(Z_r)\frac{\partial }{\partial \Theta } \log P\big (T,d \,| \, \widehat{\Lambda}(\Theta,F)  , \theta_r , \delta \big)\nonumber \\
   & = &     
    \sum_{r=1}^R\gamma(Z_r)\frac{\partial }{\partial \Theta } \left\{d  \left(\log \frac{ E dN(T) }{E  Y(T) \sum_{r=1}^{R}  \gamma(Z_{r}) \exp(\theta_r \delta_0 + X\delta_1)}+ \theta_r \delta_0 + X\delta_1 \right)\right. \nonumber\\
    & &\left. -  \exp(\theta_r \delta_0 + X\delta_1) \int_0^T  \frac{ E dN(u) }{E  Y(u) \sum_{r=1}^{R}  \gamma(Z_{r}) \exp(\theta_r \delta_0 + X\delta_1)} \right\}\nonumber \\
 & = &  \sum_{r=1}^R\gamma(Z_r) d\left[\left(\begin{array}{c}\delta_0 \\ \theta_r \\ X\end{array}\right) -\frac{M_{1}(T)}{M_{0}(T)}\right] \nonumber\\
& &   - \sum_{r=1}^R\gamma(Z_r) \left( \exp(\theta_r \delta_0 + X\delta_1)\int_0^T \left[\left(\begin{array}{c}\delta_0 \\ \theta_r \\ X\end{array}\right) 
 -\frac{M_{1}(u)}{M_{0}(u)}\right] d\Lambda(u) \right) 
    \end{eqnarray}
    where we used equation (\ref{EdN}).
The last expression is the efficient score function in the survival part of the model derived in equation (\ref{eff_score_joint}), Appendix 2.

\section*{Appendix 3: Derivation of Efficient score function in the joint Model}

In this appendix, we derive the efficient score function in the joint model using (\ref{eff_score}).
We denote $P_{r,\Theta,\Lambda}(T,d)=P\big (T,d \,| \, \Lambda  , \theta_r , \delta \big)$.

The survival part of log-likelihood function for a one observation  is
\begin{eqnarray*}
   \sum_{r=1}^R\gamma(Z_r)\log P_{r,\Theta,\Lambda}(T,d)=    \sum_{r=1}^R\gamma(Z_r)\left\{d  \big(\log \lambda(T)+ \theta_r \delta_0 + X\delta_1 \big)
   - \Lambda(T)  \exp(\theta_r \delta_0 + X\delta_1)\right\}.
\end{eqnarray*}
The score function for $ \Theta$ is
\begin{eqnarray*}
\dot{\ell}_{\Theta,\Lambda} & = &  \sum_{r=1}^R\gamma(Z_r)\frac{\partial}{\partial \Theta}   \log P_{r,\Theta,\Lambda}(T,d)
 =      \sum_{r=1}^R\gamma(Z_r)\left(\begin{array}{c}\delta_0 \\ \theta_r \\ X\end{array}\right) \left\{d    - \Lambda(T)  \exp(\theta_r \delta_0 + X\delta_1)\right\}.
\end{eqnarray*}

Let $h:[0,\tau]\rightarrow R$ be a function on $[0,\tau]$.
The path defined by 
$$d\Lambda_s=(1+sh)d\Lambda$$
is a submodel passing through $\Lambda$ at $s=0$.
The corresponding path for the $\lambda$ is
$$\lambda_s(t)=\frac{d\Lambda_s(t)}{dt}=(1+sh)\lambda(t).$$

The derivative of the log-likelihood function
\begin{eqnarray*}
   \sum_{r=1}^R\gamma(Z_r)\log P_{r,\Theta,\Lambda_s}(T,d)
   =    \sum_{r=1}^R\gamma(Z_r)\left\{d  \big(\log \lambda_s(T)+ \theta_r \delta_0 + X\delta_1 \big)
   - \Lambda_s(T)  \exp(\theta_r \delta_0 + X\delta_1)\right\}.
\end{eqnarray*}
with respect to $s$ at $s=0$ is the score operator for $\Lambda$:
$$B_{\Theta,\Lambda}h
=\frac{d}{ds}\bigg|_{s=0}\sum_{r=1}^R\gamma(Z_r)\log P_{r,\Theta,\Lambda_s}(T,d)
= \sum_{r=1}^R\gamma(Z_r) \left(d h(T)-  \exp(\theta_r \delta_0 + X\delta_1)\int_0^T h(u) d\Lambda(u) \right) .$$

\noindent {\bf Information operator $B^*_{\Theta,\Lambda}B_{\Theta,\Lambda}$}

For functions $g,h:[0,\tau]\rightarrow R$,
define a paths $d\Lambda_{s,t}=(1+sg+th+stgh)d\Lambda$.
Then
\begin{eqnarray*}
 \frac{\partial}{\partial s}\bigg|_{(s,t)=(0,0)}\sum_{r=1}^R\gamma(Z_r)\log P_{r,\Theta,\Lambda_{s,t}}(T,d) = B_{\Theta,\Lambda}g\\
 \frac{\partial}{\partial t}\bigg|_{(s,t)=(0,0)}\sum_{r=1}^R\gamma(Z_r)\log P_{r,\Theta,\Lambda_{s,t}}(T,d)= B_{\Theta,\Lambda}h.
\end{eqnarray*}
Using these we have 
\begin{eqnarray*}
\langle B_{\Theta,\Lambda}g, B_{\Theta,\Lambda}h  \rangle_{L_2(P)} 
&= &  E\{(B_{\Theta,\Lambda}g)(B_{\Theta,\Lambda}h)\}\\
 & = & -E\left\{ \frac{\partial^2}{\partial t \partial s}\bigg|_{(s,t)=(0,0)}\sum_{r=1}^R\gamma(Z_r)\log P_{r,\Theta,\Lambda_{s,t}}(T,d)\right\}\\
 & = & -E\left\{ \frac{\partial}{\partial t }\bigg|_{t=0}B_{\Theta,\Lambda_{0,t}}g\right\}\\
& = & E\left\{\sum_{r=1}^R\gamma(Z_r) \exp(\theta_r \delta_0 + X\delta_1) \int_0^{\tau}I(u \leq T) g(u)h(u) d\Lambda(u)\right\}\\ 
& = & \int_0^{\tau}h(u)  E\sum_{r=1}^R\gamma(Z_r) \exp(\theta_r \delta_0 + X\delta_1)I(u \leq T) g(u)d\Lambda(u)\\ 
& = & \left\langle E\sum_{r=1}^R\gamma(Z_r)  \exp(\theta_r \delta_0 + X\delta_1) I(u \leq T)  g(u),\ h(u)  \right\rangle_{L_2(\Lambda)}
\end{eqnarray*}

Since
$$\langle B_{\Theta,\Lambda}g, B_{\Theta,\Lambda}h  \rangle_{L_2(P)}=\langle B_{\Theta,\Lambda}^*B_{\Theta,\Lambda}g, h  \rangle_{L_2(\Lambda)},$$
we have the information operator
$$ B_{\Theta,\Lambda}^*B_{\Theta,\Lambda}g
=E \sum_{r=1}^R\gamma(Z_r)  \exp(\theta_r \delta_0 + X\delta_1) I(t \leq T)  g(t).$$
Since the operator multiplies a number, the inverse is
$$ (B_{\Theta,\Lambda}^*B_{\Theta,\Lambda})^{-1}g=\left[E\sum_{r=1}^R\gamma(Z_r)  \exp(\theta_r \delta_0 + X\delta_1) I(t \leq T)\right]^{-1} g(t).$$

\noindent {\bf Calculation of $B^*_{\Theta,\Lambda}\dot{\ell}_{\Theta,\Lambda}$}
 
 Consider a paths $(s,t) \rightarrow (\Theta+s {\bf a}, \Lambda_t)$ with $d\Lambda_t=(1+th)d\Lambda$.
 Then
 \begin{eqnarray*}
  \frac{\partial}{\partial s}\bigg|_{(s,t)=(0,0)}\sum_{r=1}^R\gamma(Z_r)\log P_{r,\Theta+s {\bf a}, \Lambda_t}(T,d) 
  & = & {\bf a}^T\dot{\ell}_{\Theta,\Lambda}\\
 \frac{\partial}{\partial t}\bigg|_{(s,t)=(0,0)}\sum_{r=1}^R\gamma(Z_r)\log P_{r,\Theta+s {\bf a}, \Lambda_t} (T,d) & = & B_{\Theta,\Lambda}h.
  \end{eqnarray*}
Using these we compute that
\begin{eqnarray*}
\langle {\bf a}^T\dot{\ell}_{\Theta,\Lambda}, B_{\Theta,\Lambda}h  \rangle_{L_2(P)} &= & 
 E\{ ({\bf a}^T\dot{\ell}_{\Theta,\Lambda})(B_{\Theta,\Lambda}h)\}\\
 & = & -E\left\{ \frac{\partial^2}{\partial t \partial s}\bigg|_{(s,t)=(0,0)}\sum_{r=1}^R\gamma(Z_r)\log P_{r,\Theta+s{\bf a},\Lambda_{t}}(T,d)\right\}\\
 & = & -E\left\{ \frac{\partial}{\partial t }\bigg|_{t=0} {\bf a}^T\dot{\ell}_{\Theta,\Lambda_t} \right\}\\
& = &  {\bf a}^TE\sum_{r=1}^R\gamma(Z_r)\left(\begin{array}{c}\delta_0 \\ \theta_r \\ X \end{array}\right) \left\{ \exp(\theta_r \delta_0 + X\delta_1) \int_0^{\tau}I(u \leq T) h(u) d\Lambda(u) \right\}\\ 
& = &  {\bf a}^T\int_0^{\tau} E\sum_{r=1}^R\gamma(Z_r)\left(\begin{array}{c}\delta_0 \\ \theta_r \\ X \end{array}\right)  \exp(\theta_r \delta_0 + X\delta_1) I(u \leq T)    h(u) d\Lambda(u)\\ 
& = &  \left\langle {\bf a}^T E\sum_{r=1}^R\gamma(Z_r)\left(\begin{array}{c}\delta_0 \\ \theta_r \\ X \end{array}\right)  \exp(\theta_r \delta_0 + X\delta_1) I(u \leq T) 
,\  h  \right\rangle_{L_2(\Lambda)}
\end{eqnarray*}

 Since
 $$\langle {\bf a}^T\dot{\ell}_{\Theta,\Lambda}, B_{\Theta,\Lambda}h  \rangle_{L_2(P)}
 = \langle {\bf a}^T B_{\Theta,\Lambda}^*\dot{\ell}_{\Theta,\Lambda}, h  \rangle_{L_2(\Lambda)},$$
we have that
$$B_{\Theta,\Lambda}^*\dot{\ell}_{\Theta,\Lambda}
= E\sum_{r=1}^R\gamma(Z_r)\left(\begin{array}{c}\delta_0 \\ \theta_r \\ X \end{array}\right)  \exp(\theta_r \delta_0 + X\delta_1) I(u \leq T) .$$

\noindent {\bf Efficient score function}

Then the efficient score function for the survival part of the model  is given by
 \begin{eqnarray}\label{eff_score_joint}
  \tilde{\ell}_{\Theta,\Lambda} 
  & = & \dot{\ell}_{\Theta,\Lambda}  - B_{\Theta,\Lambda}(B_{\Theta,\Lambda}^*B_{\Theta,\Lambda})^{-1}B_{\Theta,\Lambda}^*\dot{\ell}_{\Theta,\Lambda}\nonumber \\
 & = &  \sum_{r=1}^R\gamma(Z_r)\left(\begin{array}{c}\delta_0 \\ \theta_r \\ X\end{array}\right) \left\{d    - \Lambda(T)  \exp(\theta_r \delta_0 + X\delta_1)\right\} 
  - B_{\Theta,\Lambda}\frac{M_{1}(t)}{M_{0}(t)}\nonumber\\
 & = &  \sum_{r=1}^R\gamma(Z_r)\left(\begin{array}{c}\delta_0 \\ \theta_r \\ X\end{array}\right) \left\{d    - \Lambda(T)  \exp(\theta_r \delta_0 + X\delta_1)\right\} \nonumber\\
 & & - \sum_{r=1}^R\gamma(Z_r) \left(d \frac{M_{1}(T)}{M_{0}(T)}-  \exp(\theta_r \delta_0 + X\delta_1)\int_0^T \frac{M_{1}(u)}{M_{0}(u)} d\Lambda(u) \right)\nonumber\\
 & = &  \sum_{r=1}^R\gamma(Z_r) d\left[\left(\begin{array}{c}\delta_0 \\ \theta_r \\ X\end{array}\right) -\frac{M_{1}(T)}{M_{0}(T)}\right]\nonumber\\ 
 & & - \sum_{r=1}^R\gamma(Z_r) \left( \exp(\theta_r \delta_0 + X\delta_1)\int_0^T \left[\left(\begin{array}{c}\delta_0 \\ \theta_r \\ X\end{array}\right) 
 -\frac{M_{1}(u)}{M_{0}(u)}\right] d\Lambda(u) \right)
\end{eqnarray}
where $M_{1}(t)$ and $M_{0}(t)$ are defined in (\ref{M1M0}).

\section*{Appendix 4: Differentiability of the function $\widehat{\Lambda}(t;\Theta,F)$ given in (\ref{hat_Lambda})}


Let 
\begin{eqnarray*}
\Psi_{\Theta,F}(\Lambda)(t)
=\int_0^t  \frac{ E_{F} dN(u) }{E_{F}  Y(u) \sum_{r=1}^{R} \gamma(Z_{r}|\Theta, \Lambda)\exp(\theta_r \delta_0 + X\delta_1)}.   
\end{eqnarray*}
Then the map $\Psi_{\Theta,F}(\Lambda)$ is differentiable with respect to the parameters $(\Theta,F, \Lambda)$.
The derivatives are denoted by 
$\frac{\partial}{\partial \Theta}\Psi_{\Theta,F}(\Lambda)$, $d_{F} \Psi_{\Theta,F}(\Lambda)$ and $d_{\Lambda} \Psi_{\Theta,F}(\Lambda)$.

From (\ref{hat_Lambda}), the function $\widehat{\Lambda}(\Theta,F)$ is the solution to the operator equation:
$$\widehat{\Lambda}(\Theta,F)=\Psi_{\Theta,F}(\widehat{\Lambda}(\Theta,F)).$$
Using this the theorem below show that the function $\widehat{\Lambda}(\Theta,F)$ is differentiable in the parameter $(\Theta,F)$.

\begin{thm}\label{theorem3} (Differentiability of $\widehat{\Lambda}(\Theta,F)$)

Let $t_{max}$ be the maximum time the observed value of $T$ can get.
Suppose 
\begin{eqnarray}\label{Diff_cond}
 \left|-  \exp(\theta_r \delta_0 + X\delta_1)
    + \frac{ \sum_{g=1}^{R} \pi_g \,P\big(Y  \,|\, \theta_g , \alpha \big) P\big (T,d \,| \, \Lambda  , \theta_g , \delta \big)
  \exp(\theta_g \delta_0 + X\delta_1)}{	\sum_{g=1}^{R} \pi_g \,P\big(Y  \,|\, \theta_g , \alpha \big) P\big (T,d \,| \, \Lambda  , \theta_g , \delta \big)}\right| < \left(\int_0^{t_{max}}  d \Lambda(u)\right)^{-1}. 
\end{eqnarray}
Then the map $(\Theta,F) \rightarrow \widehat{\Lambda}(\Theta,F)$ given in (\ref{hat_Lambda})  is differentiable:
the derivatives are given by
\begin{eqnarray*}
\frac{\partial}{\partial \Theta}\widehat{\Lambda}(\Theta,F)& =  &[I-d_{\Lambda} \Psi_{\Theta,F}(\Lambda) ]^{-1}\frac{\partial}{\partial \Theta}\Psi_{\Theta,F}(\Lambda),\\
d_F \widehat{\Lambda}(\Theta,F)& =  &[I-d_{\Lambda} \Psi_{\Theta,F}(\Lambda) ]^{-1}d_F\Psi_{\Theta,F}(\Lambda).
\end{eqnarray*}

\end{thm}

{\bf Proof.}
For each fixed $t$,  the maps $\Theta \rightarrow \Psi_{\Theta,F}(\Lambda)(t)$ and 
$\Lambda \rightarrow \Psi_{\Theta,F}(\Lambda)(t)$  are differentiable real valued maps,
by the result in Appendix 4, there are some $\Theta^*$ and $\Lambda^*$ with $\|\Theta^*-\Theta\| \leq \|\Theta'- \Theta\|$ 
and
$\|\Lambda^*-\widehat{\Lambda}(\Theta,F)\| \leq \|\widehat{\Lambda}(\Theta',F)-\widehat{\Lambda}(\Theta,F)\|$ 
such that
\begin{eqnarray*}
& & \widehat{\Lambda}(\Theta',F)-\widehat{\Lambda}(\Theta,F)\\
& = & \Psi_{\Theta',F}(\widehat{\Lambda}(\Theta',F))- \Psi_{\Theta,F}(\widehat{\Lambda}(\Theta,F)) \\
& = & \Psi_{\Theta',F}(\widehat{\Lambda}(\Theta',F))- \Psi_{\Theta,F}(\widehat{\Lambda}(\Theta',F)) 
  + \Psi_{\Theta,F}(\widehat{\Lambda}(\Theta',F))- \Psi_{\Theta,F}(\widehat{\Lambda}(\Theta,F)) \\
& = & \frac{\partial}{\partial \Theta} \Psi_{\Theta^*,F}(\widehat{\Lambda}(\Theta',F)) (\Theta'- \Theta)
  +d_{\Lambda}  \Psi_{\Theta,F}(\Lambda^*) (\widehat{\Lambda}(\Theta',F)-\widehat{\Lambda}(\Theta,F)).
\end{eqnarray*}
It follows that
\begin{eqnarray*}
[I-d_{\Lambda} \Psi_{\Theta,F}(\Lambda^*) ](\widehat{\Lambda}(\Theta',F)-\widehat{\Lambda}(\Theta,F) )
& = &\frac{\partial}{\partial \Theta}\Psi_{\Theta^*,F}(\widehat{\Lambda}(\Theta',F))(\Theta'-\Theta).
\end{eqnarray*}

If we can show the inverse $[I-d_{\Lambda} \Psi_{\Theta,F}(\Lambda^*) ]^{-1}$ exists, then
\begin{eqnarray*}
\widehat{\Lambda}(\Theta',F)-\widehat{\Lambda}(\Theta,F) 
& = &[I-d_{\Lambda} \Psi_{\Theta,F}(\Lambda^*) ]^{-1}\frac{\partial}{\partial \Theta}\Psi_{\Theta^*,F}(\widehat{\Lambda}(\Theta',F))(\Theta'-\Theta)\\
& = &[I-d_{\Lambda} \Psi_{\Theta,F}(\Lambda) ]^{-1}\frac{\partial}{\partial \Theta}\Psi_{\Theta,F}(\widehat{\Lambda}(\Theta,F))(\Theta'-\Theta)+o(\Theta'-\Theta),
\end{eqnarray*}
as $\Theta' \rightarrow \Theta$.
Thus the derivative of $\widehat{\Lambda}(\Theta,F)$ with respect to $\Theta$ is
$$\frac{\partial}{\partial \Theta}\widehat{\Lambda}(\Theta,F)=[I-d_{\Lambda} \Psi_{\Theta,F}(\Lambda) ]^{-1}\frac{\partial}{\partial \Theta}\Psi_{\Theta,F}(\Lambda).$$

Now we show the inverse $[I-d_{\Lambda} \Psi_{\Theta,F}(\Lambda) ]^{-1}$ exists.
This follows if we show $$\|d_{\Lambda} \Psi_{\Theta,F}(\Lambda)\|<1$$ in the operator norm.

Using (\ref{EdN}), at the true value of parameters, the derivative $d_{\Lambda} \Psi_{\Theta,F}(\Lambda)(t)$ can be expressed as 
\begin{eqnarray*}
d_{\Lambda} \Psi_{\Theta,F}(\Lambda)(t)\{\Lambda'-\Lambda\} =
 -\int_0^t  d \Lambda(u) 
 \frac{  E_{F}  Y(u) \sum_{r=1}^{R} \gamma(Z_{r}|\Theta, \Lambda)
 \exp(\theta_r \delta_0 + X\delta_1) \frac{d_{\Lambda}\gamma(Z_{r}|\Theta, \Lambda) (\Lambda'-\Lambda) }{\gamma(Z_{r}|\Theta, \Lambda)}  }{E_{F}  Y(u) \sum_{r=1}^{R} \gamma(Z_{r}|\Theta, \Lambda)\exp(\theta_r \delta_0 + X\delta_1)}
\end{eqnarray*}
where
\begin{eqnarray*}
& & \frac{d_{\Lambda} \gamma(Z_{r}|\Theta, \Lambda)(\Lambda'-\Lambda)}{\gamma(Z_{r}|\Theta, \Lambda)} \\
&  = & 
   \left[-  \exp(\theta_r \delta_0 + X\delta_1)
    + \frac{ \sum_{g=1}^{R} \pi_g \,P\big(Y  \,|\, \theta_g , \alpha \big) P\big (T,d \,| \, \Lambda  , \theta_g , \delta \big)
  \exp(\theta_g \delta_0 + X\delta_1)}{	\sum_{g=1}^{R} \pi_g \,P\big(Y  \,|\, \theta_g , \alpha \big) P\big (T,d \,| \, \Lambda  , \theta_g , \delta \big)}\right] \{\Lambda'(T)-\Lambda(T)\}.
\end{eqnarray*}

By the assumption (\ref{Diff_cond}), we have that
\begin{eqnarray*}
|d_{\Lambda} \Psi_{\Theta,F}(\Lambda)(t)\{\Lambda'-\Lambda\}| \leq 
   \left(\int_0^{t_{max}}  d \Lambda(u)\right)  \left(\int_0^{t_{max}}  d \Lambda(u)\right)^{-1}  \|\Lambda'(T)-\Lambda(T)\|
\leq   \|\Lambda'(T)-\Lambda(T)\|.
\end{eqnarray*}
This shows $\|d_{\Lambda} \Psi_{\Theta,F}(\Lambda)\|<1$ in the operator norm at the true value of the parameters. We assume this holds in some neighborhood of the true values.
It follows that, in the neighborhood, the function $\widehat{\Lambda}(\Theta,F)$ is differentiable with respect to $\Theta$.

A similar proof can show the differentiability of the function $\widehat{\Lambda}(\Theta,F)$ with respect to the parameter $F$.

\section*{Appendix 5: Mean value theorem for functional}

Suppose $f: H \rightarrow  R$ is a Hadamard differentiable  real valued map on a convex subset $H$ of a Banach space.
We denote the derivative by $df$.
Suppose $\eta ,\eta_0 \in H$, then there exist $\eta^* \in H$ such that $\|\eta^* - \eta_0 \| \leq  \|\eta - \eta_0 \|$ and
$$f(\eta) = f(\eta_0)+df(\eta^*)(\eta-\eta_0).$$

Moreover if $\eta$ is partitioned $\eta=(\eta_1,\eta_2)$ so that $\eta_0=(\eta_{01},\eta_{02})$ and
$df=(df_1,df_2)$ with $df_i$ the derivative with respect to $\eta_i$, $i=1,2$.

Then there are $\eta_1^*$ and $\eta_2^*$ such that
$\|\eta^*_1 - \eta_{01} \| \leq  \|\eta_1 - \eta_{01} \|$ and $\|\eta^*_2 - \eta_{02} \| \leq  \|\eta_2 - \eta_{02} \|$, and
$$f(\eta_1,\eta_2) = f(\eta_{01},\eta_{02})+df_1(\eta^*_1,\eta^*_2)(\eta_1-\eta_{01})+df_2(\eta^*_1,\eta^*_2)(\eta_2-\eta_{02}).$$

{\bf Proof}
Since $H$ is a convex set $\eta_0 + t(\eta-\eta_0) \in H$ for each $t \in [0,1]$.   
Then the map $t \rightarrow f(\eta_0 + t(\eta-\eta_0))$ is a differentiable real valued function. By the mean value theorem there is a $t^* \in (0,1)$ such that
$$f(\eta) = f(\eta_0)+df(\eta_0 + t^*(\eta-\eta_0)) (\eta-\eta_0).$$ 
Let $\eta^*=\eta_0 + t^*(\eta-\eta_0)$ and, for the partitioned version this is 
$(\eta^*_1,\eta^*_2)=(\eta_{01},\eta_{02})+t^*[(\eta_{1},\eta_{2})-(\eta_{01},\eta_{02})]$, then the desired result follows.

\end{document}